\documentclass[12pt,twoside,english]{article}
\usepackage{amsfonts,babel,amssymb,amsmath,amsthm}
\usepackage{pstricks-add}
\usepackage{slashbox}

\numberwithin{equation}{section}

\newtheorem{proposition}{Proposition}[section]

\newtheorem{theorem}[proposition]{Theorem}

\newtheorem{remark}{Remark}[section]

\def\div{\mathop{\mathrm{div}}\nolimits}

\title{A decoupled preconditioning technique for 
a mixed Stokes-Darcy model\footnote{This research was partially supported 
by  Ministery of Education of Spain through the Project MTM2010-18427. }}

\author{
{\sc Antonio M\'arquez}\thanks{Departamento de Construcci\'on e Ingenier\'\i a de
Fabricaci\'on, Universidad de Oviedo, Oviedo, Espa\~na,
e-mail: {\tt amarquez@uniovi.es}},
$\,\,$
{\sc Salim Meddahi}\thanks{Departamento de Matem\'aticas, Facultad de Ciencias,
Universidad de Oviedo, Calvo Sotelo s/n, Oviedo, Espa\~na,
e-mail: {\tt salim@uniovi.es}}, \\ and \\
{\sc Francisco-Javier Sayas}\thanks{Department of Mathematical Sciences,
University of Delaware, Newark DE 19716, USA, e-mail: {\tt
fjsayas@math.udel.edu}}
}

\date{}

\begin{document}

\maketitle

\begin{abstract}
We propose an efficient iterative method to solve 
the mixed Stokes-Dracy model for coupling fluid and porous media flow.  
The weak formulation of this problem leads to a coupled, 
indefinite, ill-conditioned and symmetric linear system of equations. We apply a  
decoupled preconditioning technique requiring only good solvers for  
the local mixed-Darcy and Stokes subproblems.
We prove that the method is asymptotically optimal and confirm, with   
numerical experiments, that the performance of the preconditioners  
does not deteriorate on arbitrarily fine meshes. 
\end{abstract}



\section{Introduction}\label{sec:1}

The Stokes-Darcy problem describes filtration processes that find many
important applications in porous media problems. 
Usually, a surface free flow of a liquid is modeled by Stokes equations and 
the flow confined in the porous media is governed by Darcy equations.  
 The interaction of the local models is commonly handled through 
 the Beavers-Joseph-Saffman (BJS) interface conditions, cf. \cite{beavers, saffman, jager}.

Recently, there has been active research on the mathematical and numerical analysis for this model.
The weak formulation of this problem is generally obtained by coupling 
the usual velocity-pressure mixed formulation in the Stokes domain with either the primal formulation 
($H^1$-approach, \cite{miglio}) or the 
mixed formulation ($\mathbf{H}(\text{div})$-approach,  \cite{layton, GatOyaSay11} ) 
in the Darcy domain. In both cases, 
optimal iterative methods are of crucial importance to solve efficiently the discrete Stokes-Darcy model
since the corresponding linear systems of algebraic equations are indefinite and
ill-conditioned. Several optimal iterative solvers 
such as the Dirichlet-Neumann or Robin-Robin 
domain decomposition methods \cite{CaoGunzburger, chen, Discacciati1, Discacciati2, Discacciati3}, 
multi-grid methods \cite{CaiMu, CaiXu, MuXu} have been proposed for the model based on an $H^1$-approach in 
the Darcy domain. The common feature for most of these methods 
consists in decoupling the global model in such a way that, only 
independent Stokes and Darcy subproblems are involved in the iterative process. 
To the authors' knowledge, the only solution procedure for 
the $\mathbf{H}(\text{div})$-conforming Darcy flux approach 
was proposed in \cite{arbogast}. In this paper, 
the problem is written as a global saddle point problem and a solver is implemented for the scheme using 
an inexact Uzawa technique relying on an expensive preconditioner. 
Our purpose here is to devise a decoupled preconditioning strategy that allows to apply  
existing optimized solvers to each local model independently. 

Recently,  a systematic way  to obtain convergent finite
element schemes for the Darcy-Stokes flow problem by combining well-known mixed
finite elements that are separately convergent for the $\mathbf{H}(\div)$-Darcy formulation and 
the Stokes problem was proposed in \cite{anteriorPaper}. We take advantage here of a 
fluid-to-pressure (FtP) operator to reinterpret in this formulation the Darcy system 
as a nonlocal boundary condition for the Stokes problem.
The corresponding discrete equations are written in terms of a symmetric and indefinite linear 
operator that enjoys the same spectral properties of the local discrete Stokes problem. 

Many different iterative methods for solving the saddle point problems that result from 
the finite element discretization of the Stokes equations are known.  
There are, for instance,  many variants of the so-called inexact Uzawa methods. 
Block triangular preconditioners for saddle point operators has also been discussed by
many authors. An example of such preconditioners has been introduced in \cite{bramblePasciak} 
by Bramble and Pasciak. In this strategy the original saddle point system is premultiplied   
by a block triangular operator and then, the resulting positive definite system is solved 
by a preconditioned conjugate gradient method. A possible difficulty is that this approach 
requires a proper choice of a critical scaling parameter to obtain a positive definite operator. 
Finally, we mention the block diagonal preconditioners for the minimal residual method (MINRES), cf.  
\cite{Elman, mardalWinther} and the references therein.  
 A comparative study of representing methods from each of these three  
classes is considered in \cite{petersReusken}. The inexact Uzawa method is not feasible 
in our case because of the nonlocal character of the discrete flux-to-pressure operator 
$C^h_{\mathrm{S}}$ appearing in the principal block of our saddle point problem 
\eqref{matrix1}. The conclusion in  \cite{petersReusken} is that 
the preconditioned MINRES method may be slower than 
the Bramble-Pasciak method but it is more robust (it even converges without preconditioning) and 
it is parameter free. 

Applying a preconditioned MINRES method in our case requires, at each iteration 
step, the solution of two local problems: a vector Laplace equation in the fluid and the mixed formulation 
of the Darcy problem in the porous media. This saddle point problem in the Darcy domain is again 
solved with a preconditioned MINRES method. 
The global algorithm has then the structure of an outer-inner MINRES iteration 
process. Thus, for this method a stopping criterion (tolerance parameter) for 
the inner iteration is needed. The preconditioner for the inner MINRES may be constructed 
by using techniques from \cite{arnoldFalkWinther, HiptmairXu}. We use here 
the nodal auxiliary space preconditioning 
technique introduced in \cite{HiptmairXu} by Hiptmair and Xu 
to solve $\mathbf{H}(\text{div})$-elliptic problems. With this choice, our decoupled iterative 
process consists in two nested MINRES methods whose preconditioners only require 
the solution of several second-order $H^1$-elliptic 
problems  in the Stokes and the Darcy domains. Standard multigrid techniques or 
domain decomposition methods can then be applied to reduce the computational effort.
Theoretical analysis and numerical experiments show the optimality  
and efficiency of the proposed decoupled iterative solver. 

The rest of the paper is organized as follows. In Section \ref{sec:2} we summarize the results of \cite{anteriorPaper} 
introducing the model problem, the variational formulation, general conditions for convergence of a Galerkin discretization 
and examples of spaces leading to convergent methods. A reinterpretation of the continuous formulation 
in terms of a Fluid-to-Pressure operator and the derivation of its discrete counterpart  
are presented in Section \ref{sec:5}. We take advantage of the equivalent formulation of the 
discrete Stokes-Darcy problem to deduce, in Section \ref{sec:6}, 
a decoupled iterative solver based on a preconditioned MINRES method.
Finally, numerical experiments are reported in Section \ref{sec:7}.

\paragraph{Notation and background.}

Boldface fonts will be used to denote vectors  and vector valued functions. Also, if $H$ 
is a vector space of scalar functions, 
$\mathbf H$ will denote the space of $\mathbb R^d$ valued functions whose components are in $H$, 
endowed with the product norm. 

Given an integer $m\geq 1$ and a bounded Lipschitz domain $\mathcal O\subset \mathbb R^d$, $(d=2,3)$,  
we denote by $\|\cdot \|_{H^m(\mathcal O)}$ the norm in the usual Sobolev space 
$H^m(\mathcal O)$, cf. \cite{AdaFou03}. For economy of notation, $(\cdot,\cdot)_{\mathcal O}$ stands for the 
inner product in $L^2(\mathcal{O})$ and $\|\cdot \|_{\mathcal{O}}$ is the corresponding norm. 
We recall that $H^{1/2}(\partial \mathcal{O})$ represents the image of $H^1(\mathcal{O})$ by the trace operator. 
Its dual with respect to the pivot space $L^2(\partial \mathcal{O})$ is denoted $H^{-1/2}(\partial \mathcal{O})$. 
For definition and basic properties of the spaces 
$\mathbf H(\mathrm{div},\mathcal O)$ and $\mathbf H(\mathbf{curl},\mathcal O)$, 
we refer to \cite{GirRav86}. We will denote by $\mathbf H_{0}(\mathrm{div}, \mathcal O)$ 
the subspace of fields from $\mathbf H(\mathrm{div}, \mathcal O)$
with zero normal trace on the boundary $\partial \mathcal O$ and by $\mathbf H_{0}(\mathbf{curl}, \mathcal O)$ 
the subspace of fields from $\mathbf H(\mathbf{curl}, \mathcal O)$
with zero tangential trace on $\partial \mathcal O$.

For $k\ge 0$, $\mathbb P_k(\mathcal O)$ will denote 
the space of $d-$variate polynomials 
of degree not greater than $k$ defined on the set $\mathcal O\subset \mathbb R^d$ with non-trivial interior.
Finally, at the discrete level, the letter $h$ (with or without geometric meaning) will be used to 
denote discretization. The expression 
$a\lesssim b$ will be used to mean that there exists $C>0$ independent of $h$ such that $a\le C\,b$ for 
all $h$. Similarly, we write $a\simeq b$ when there exist constants $C>c>0$ independent of $h$ such that 
$c a \leq b \leq C a$.

Let us consider a linear operator $A^h:\: X_h \to X_h^*$ acting between a  
finite dimensional subspace $X_h$ of a Hilbert space $X\subset L^2(\mathcal O)$ and its dual $X_h^*$. 
Assume that we have chosen a basis in $X_h$ and  that the coefficients 
of $u_h, v_h \in X_h$  in this basis are given by $\bar u, \bar v\in \mathbb{R}^n$, 
where $n$ is the dimension of $X_h$. 
We define the matrix realization $\mathbf{A}^h\in \mathbb{R}^{n \times n}$ of $A^h$ by
\begin{equation}\label{realization}
 \langle\mathbf{A}^h \bar u , \bar v\rangle_2 =  \langle A^h  u_h , v_h\rangle_{X_h^*\times X_h} 
 \quad \forall u_h, v_h \in X_h,
\end{equation}
where $\langle\cdot , \cdot \rangle_2$ stands for the Euclidean scalar product in $\mathbb{R}^n$. 
Moreover, if $I^h:\: X_h \to X_h'$ is the Riesz operator given by
\[
 \langle I^h u_h ,  v_h \rangle_{X_h^*\times X_h} =  (  u_h , v_h)_{\mathcal O} 
 \quad \forall u_h, v_h \in X_h,
\]
then, the corresponding matrix realization $\mathbf{M}^h\in \mathbb{R}^{n \times n}$, 
usually  referred to as the mass matrix, is defined by
\[
 \langle\mathbf{M}^h \bar u , \bar v\rangle_2 =  \langle I^h  u_h , v_h\rangle_{X_h^*\times X_h} 
 \quad \forall u_h, v_h \in X_h.
\]

\section{Statement of the problem and discretization}\label{sec:2}

Let us consider a domain $\Omega\subset \mathbb
R^d$ ($d=2$ or $d=3$) with polyhedral Lipschitz boundary. We assume that $\Omega$ is
subdivided into two subdomains by a Lipschitz polyhedral interface
$\Sigma$. The subdomains are denoted $\Omega_{\mathrm S}$ and
$\Omega_{\mathrm D}$ (S stands for Stokes and D for Darcy). We also
denote $\Gamma_{\mathrm S}:=\partial\Omega_{\mathrm S}\setminus\Sigma$ and $\Gamma_{\mathrm D}:=\partial\Omega_{\mathrm D}
\setminus\Sigma.$
The normal vector field $\mathbf{n}$ on $\partial\Omega$ is
chosen to point outwards. We also denote by $\mathbf{n}$ the normal vector on $\Sigma$ 
that points from $\Omega_{\mathrm S}$ to $\Omega_{\mathrm D}$. 

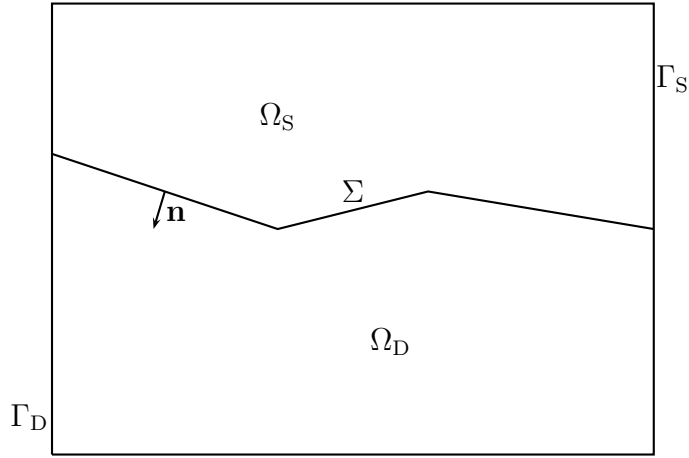
\begin{center}
 \begin{pspicture}(-4,-3)(4,3)
 \psline(-4,-3)(-4,3)(4,3)(4,-3)(-4,-3)
  \psline(-4,1)(-1,0)(1,0.5)(4,0)
  \psline{->}(-2.5,0.5)(-2.65,0)
  \rput(-2.35,0.2){$\mathbf{n}$}
  \rput(0.5,-1.5){$\Omega_{\mathrm D}$}
  \rput(-1,1.5){$\Omega_{\mathrm S}$}
  \rput(0,0.5){$\Sigma$}
  \rput(4.25,2){$\Gamma_{\mathrm S}$}
  \rput(-4.3,-2.5){$\Gamma_{\mathrm D}$}
 \end{pspicture}
\end{center}

\subsection{Variational formulation}
In the region $\Omega_{\mathrm S}$, 
the fluid flow is assumed to satisfy the Stokes system 
\begin{equation}\label{eq:00.1}
-\textbf{div} \left( 2\nu \boldsymbol\varepsilon(\mathbf u_{\mathrm S})  -p_{\mathrm S} \mathbf{I} \right)  = \mathbf f_{\mathrm S}, 
\qquad \mathrm{div}\,\mathbf u_{\mathrm S} = 0 
\qquad \mbox{in $\Omega_{\mathrm S}$},
\end{equation}
where $\mathbf{I}$ is  the identity in $\mathbb{R}^d$ and 
$\boldsymbol\varepsilon(\mathbf u_{\mathrm S}):=\frac12 (\boldsymbol\nabla\mathbf 
u_{\mathrm S} + \boldsymbol\nabla^\top\mathbf u_{\mathrm S})$ is the deformation tensor, 
$\nu>0$ is the kinematic viscosity and $\mathbf f_{\mathrm S}$ is the external body force. 
In the porous region $\Omega_{\mathrm D}$, the governing equations are given by the following Darcy system 
\begin{equation}\label{eq:00.2}
\mathbf K^{-1}\mathbf u_{\mathrm D}+\nabla p_{\mathrm D}=\mathbf 0, \qquad \mathrm{div}\,\mathbf u_{\mathrm D}= f_{\mathrm D}
\qquad \mbox{in $\Omega_{\mathrm D}$},
\end{equation}
where $f_{\mathrm D}$ is the source (or sink) term and the hydraulic conductivity tensor of the porous medium 
$\mathbf K(\mathbf x)$ is symmetric and uniformly bounded and positive definite, i.e., 
\[
 0< k_1 | \boldsymbol{\xi} |^2 \leq \mathbf K(\mathbf x)\boldsymbol{\xi}\cdot \boldsymbol{\xi} 
 \leq k_2 | \boldsymbol{\xi} |^2 \quad \text{for a.e. $\mathbf{x}\in \Omega_{\mathrm D}$}, \quad \forall 
 \boldsymbol{\xi}\in \mathbb{R}^d,
\]
for some constants $k_2 \geq k_1>0$. 
On the outer boundaries we consider the 
homogeneous (essential) boundary conditions
\begin{equation}\label{eq:00.3}
\mathbf u_{\mathrm S}=\mathbf 0 \quad\mbox{on $\Gamma_{\mathrm S}$}, \qquad \mbox{and}\qquad \mathbf 
u_{\mathrm D}\cdot\mathbf{n} = 0 \quad \mbox{on $\Gamma_{\mathrm D}$},
\end{equation}
and on the interface between the fluid and porous media regions we impose conditions ensuring 
mass conservation, balance of normal forces and the Beavers-Joseph-Saffman condition \cite{beavers, saffman},
\begin{equation}\label{eq:00.4}
\mathbf u_{\mathrm S}\cdot\mathbf{n} = \mathbf u_{\mathrm D}\cdot\mathbf{n}, \qquad 2\nu 
\boldsymbol\varepsilon(\mathbf u_{\mathrm S})\mathbf{n} - p_{\mathrm S}\mathbf{n} + \kappa \boldsymbol\pi_t 
\mathbf u_{\mathrm S} = -p_{\mathrm D}\mathbf{n} \qquad\mbox{on $\Sigma$},
\end{equation}
where $\boldsymbol\pi_t \mathbf w:=\mathbf w-(\mathbf w\cdot\mathbf{n})\mathbf{n}$ 
and $\kappa$ is a positive and bounded  function  depending on $\mathbf K$, $\nu$, and an experimentally determined 
friction constant, cf. \cite{beavers, CaiMu, CaoGunzburger}. 
 
Because of the mass conservation 
condition across $\Sigma$, the homogeneous Dirichlet boundary condition for $\mathbf u_{\mathrm S}$ on $\Gamma_{\mathrm S}$ 
and the incompressibility condition in the Stokes domain, we can easily show that
\begin{equation}\label{eq:0}
\int_{\Omega_{\mathrm D}} f_{\mathrm D}=0
\end{equation}
is a necessary condition for existence of solution. The pressure field is defined up to an additive constant. 
We will normalize it by imposing that 
\[
\int_{\Omega_{\mathrm D}} p_{\mathrm D}=0.
\]
For the velocity field, we will use the space
\[
\mathbb X:=\{ \mathbf u=(\mathbf u_{\mathrm S},\mathbf u_{\mathrm
D})\in \mathbf H^1_{\mathrm S}(\Omega_{\mathrm S}) \times \mathbf
H_{\mathrm D}(\mathrm{div},\Omega_{\mathrm D})\,:\, \mathbf
u_{\mathrm S}\cdot\mathbf{n}=\mathbf u_{\mathrm
D}\cdot\mathbf{n} \quad \mbox{on $\Sigma$}\} \subset \mathbf
H_0(\mathrm{div},\Omega),
\]
where  
\begin{eqnarray}\label{eq:1.10}
\mathbf H^1_{\mathrm S}(\Omega_{\mathrm S}) &:=& \{ \mathbf u \in
\mathbf{H}^1(\Omega_{\mathrm{S}}) \,:\, \mathbf u=\mathbf 0 \quad \mbox{on
$\Gamma_{\mathrm S}$}\}\\
\label{eq:1.11} \mathbf H_{\mathrm D}(\mathrm{div}, \Omega_{\mathrm
D}) &:=& \{ \mathbf u \in \mathbf H(\mathrm{div},\Omega_{\mathrm
D})\,:\, \mathbf u \cdot \mathbf{n} =0 \quad \mbox{on
$\Gamma_{\mathrm D}$}\}.
\end{eqnarray}
The space $\mathbb X$ will be endowed with the product norm.
The space for the pressure field is $\mathbb Q:= L^2(\Omega_{\mathrm S}) \times L^2_0(\Omega_{\mathrm D}),$
where $L^2_0(\mathcal O):= \{ p \in L^2(\mathcal O)\,:\, (p,1)_{\mathcal
O}=0\}.$
The pressure field is represented as $p:=(p_{\mathrm S}, p_{\mathrm D})\in \mathbb{Q}$. 
Adding an appropriate constant in a postprocessing step, the normalization condition $(p,1)_{\Omega_{\mathrm D}}=0$ 
can be modified to $(p,1)_\Omega=0$.
The space $\mathbb Q$ is endowed with the corresponding
product norm.

We consider four bilinear forms, two in the Stokes domain and two in the Darcy domain:
\begin{eqnarray}\label{eq:aS}
 a_{\mathrm{S}}(\mathbf{u}_{\mathrm{S}}, \mathbf{u}_{\mathrm{S}}) &:=& 2\nu(\boldsymbol\varepsilon(\mathbf
u_{\mathrm S}),\boldsymbol\varepsilon(\mathbf v_{\mathrm
S}))_{\Omega_{\mathrm S}}+ \langle\kappa
\boldsymbol\pi_t\mathbf u_{\mathrm S},\boldsymbol\pi_t\mathbf
v_{\mathrm S}\rangle_\Sigma,\\
 a_{\mathrm{D}}(\mathbf{u}_{\mathrm{D}}, \mathbf{u}_{\mathrm{D}}) &:=& 
 (\mathbf K^{-1} \mathbf u_{\mathrm
D},\mathbf v_{\mathrm D})_{\Omega_{\mathrm D}},\\
 b_{\mathrm{S}}(\mathbf u_{\mathrm{S}},q_{\mathrm{S}}) &:=& 
 (\mathrm{div}\,\mathbf u_{\mathrm{S}},q_{\mathrm{S}})_{\Omega_{\mathrm{S}}},\\
 b_{\mathrm{D}}(\mathbf u_{\mathrm{D}},q_{\mathrm{D}}) &:=& 
 (\mathrm{div}\,\mathbf u_{\mathrm{D}},q_{\mathrm{D}})_{\Omega_{\mathrm{D}}}.
\end{eqnarray}
 These bilinear forms are combined to build the 
diagonal bilinear form of the mixed problem $a: \mathbb X \times \mathbb X\to \mathbb R$, given by
\[
a(\mathbf u,\mathbf v):= a_{\mathrm{S}}(\mathbf{u}_{\mathrm{S}}, \mathbf{u}_{\mathrm{S}})+
a_{\mathrm{D}}(\mathbf{u}_{\mathrm{D}}, \mathbf{u}_{\mathrm{D}}),
\]
as well as  $b:\mathbb X \times \mathbb Q \to \mathbb
R$ given by
\begin{equation}\label{eq:1.0a}
b(\mathbf u,q):= b_{\mathrm{S}}(\mathbf u_{\mathrm{S}},q_{\mathrm{S}}) + 
b_{\mathrm{D}}(\mathbf u_{\mathrm{D}},q_{\mathrm{D}}).
\end{equation}
A well posed variational form of the Darcy-Stokes problem (cf. \cite[Proposition 2.3]{anteriorPaper})
is: find $(\mathbf u,p) \in \mathbb X \times \mathbb Q$
such that
\begin{equation}\label{eq:1.vp}
\begin{array}{rll}
a(\mathbf u,\mathbf v) - b(\mathbf v,p) &=(\mathbf
f_{\mathrm S},\mathbf v_{\mathrm S})_{\Omega_{\mathrm S}} & \forall
\mathbf v \in \mathbb X,\\[1.5ex]
b(\mathbf u,q) &=(f_{\mathrm D},q_{\mathrm
D})_{\Omega_{\mathrm D}} & \forall q\in \mathbb Q.
\end{array}
\end{equation}

\subsection{The discrete problem}\label{sec:3}

We start by creating shape-regular triangulations $\{\mathcal T^h_{\mathrm S}\}_h$ and $\{\mathcal
T^h_{\mathrm D}\}_h$  of $\overline{\Omega}_{\mathrm S}$ and $\overline{\Omega}_{\mathrm D}$
respectively, consisting of triangles (tetrahedra in the three dimensional case) of diameter not larger than $h$.
The triangulations create two inherited
partitions of $\Sigma$, respectively denoted $\Sigma^h_{\mathrm S}$
and $\Sigma^h_{\mathrm D}$. 
Let us consider finite dimensional subspaces of piecewise polynomial) to approximate velocity
and pressure in the Stokes domain
\[
\mathbf V^h(\Omega_{\mathrm S}) \subset \mathbf H^1(\Omega_{\mathrm
S}), \qquad  L^h_0(\Omega_{\mathrm S})\subset L^2_0(\Omega_{\mathrm
S}),\qquad L^h(\Omega_{\mathrm S})=L^h_0(\Omega_{\mathrm S})\oplus
\mathbb P_0(\Omega_{\mathrm S}),
\]
as well as in the Darcy domain
\[
\mathbf H^h(\Omega_{\mathrm D}) \subset \mathbf
H(\mathrm{div},\Omega_{\mathrm D}), \qquad  L^h_0(\Omega_{\mathrm
D})\subset L^2_0(\Omega_{\mathrm D}),\qquad L^h(\Omega_{\mathrm
D})=L^h_0(\Omega_{\mathrm D})\oplus \mathbb P_0(\Omega_{\mathrm D}).
\]
We also need to consider the spaces
\begin{equation}\label{eq:3.4}
\mathbf V^h_{\mathrm S}(\Omega_{\mathrm S}):= \mathbf
H^h(\Omega_{\mathrm S}) \cap \mathbf H^1_{\mathrm S}(\Omega_{\mathrm
S}), \qquad \mathbf H^h_{\mathrm D}(\Omega_{\mathrm D}):= \mathbf
H^h(\Omega_{\mathrm D}) \cap \mathbf H_{\mathrm
D}(\mathrm{div},\Omega_{\mathrm D}),
\end{equation}
as well as the discrete spaces of normal components on $\Sigma$, namely, 
\begin{eqnarray*}
\Phi^h_{\mathrm S}&:=& \{ \mathbf u_h\cdot \mathbf{n} \,:\, \mathbf u_h
\in \mathbf V^h_{\mathrm S}(\Omega_{\mathrm S})\}\subset L^2(\Sigma),\\
\Phi^h_{\mathrm D}&:=& \{ \mathbf u_h\cdot \mathbf{n} \,:\, \mathbf u_h
\in \mathbf H^h_{\mathrm D}(\Omega_{\mathrm D})\}\subset L^2(\Sigma).
\end{eqnarray*}
We will assume that $\Phi^h_{\mathrm D}$ contains at least the space of piecewise constant 
functions on $\Sigma^h_{\mathrm{D}}$ and denote by $R^h_D$ the $L^2(\Sigma)$-orthogonal 
projection onto $\Phi^h_{\mathrm D}$.

The method we are proposing is a Galerkin discretization of the
variational problem \eqref{eq:1.vp} using the spaces
\begin{eqnarray*}
\mathbb X^h &:=&  \{ \mathbf u_h \equiv (\mathbf
u^h_{\mathrm S},\mathbf u^h_{\mathrm D}) \in \mathbf V^h_{\mathrm
S}(\Omega_{\mathrm S})\times \mathbf H^h_{\mathrm D}(\Omega_{\mathrm
D})\,:\,  \mathbf u^h_{\mathrm D}\cdot\mathbf{n} = R^h_D(\mathbf u^h_{\mathrm S}\cdot\mathbf{n})\,\, 
\text{on $\Sigma$}\},\\
\mathbb Q^h &:=& L^h(\Omega_{\mathrm S})\times L^h_0(\Omega_{\mathrm D}),
\end{eqnarray*}
that is, we look for $(\mathbf u_h,p_h)\in \mathbb X^h
\times \mathbb Q^h$ such that
\begin{equation}\label{eq:2.0}
\begin{array}{rll}
a(\mathbf u_h,\mathbf v_h)-b(\mathbf v_h,p_h)&=(\mathbf
f_{\mathrm S},\mathbf v_h)_{\Omega_{\mathrm S}} & \forall \mathbf
v_h \in \mathbb X^h,\\[1.5ex]
b(\mathbf u_h,q_h)&=(f_{\mathrm D}, q_h)_{\Omega_{\mathrm
D}} & \forall q_h \in \mathbb Q^h.
\end{array}
\end{equation}

\begin{remark}
Note that  $\mathbb X^h \not\subset \mathbb X$ unless $\Phi^h_{\mathrm S}\subset \Phi^h_{\mathrm D}$ 
in which case \eqref{eq:2.0} becomes a conforming Galerkin approximation  of \eqref{eq:1.vp}. We will say 
that the discretization is conforming if $\Phi^h_{\mathrm S}\subset \Phi^h_{\mathrm D}$ 
(and therefore $\mathbb X^h\subset \mathbb X$) and non-conforming otherwise. 
\end{remark}
 
The following result is proved in \cite[Proposition 3.2]{anteriorPaper}. Inf-sup conditions are 
written in terms of the spaces
\[
\mathbf V^h_0(\Omega_{\mathrm S}):= \mathbf V^h(\Omega_{\mathrm
S})\cap \mathbf H^1_0(\Omega_{\mathrm S}), \qquad \mathbf
H^h_0(\Omega_{\mathrm D}):=\mathbf H^h(\Omega_{\mathrm D}) \cap
\mathbf H_0(\mathrm{div},\Omega_{\mathrm D}),
\]
which arise from the application of
the discretization method to problems with homogeneous boundary
conditions on the entire boundary of each subdomain.

\begin{theorem}\label{thm:main}
Let us assume that there exist a linear operator $\mathbf L_h:\Phi_{\mathrm D}\to \mathbf H^h_{\mathrm D}(\Omega_{\mathrm D})$ 
and a general positive constant $\beta$, independent of $h$, such that:
\begin{equation}\label{eq:2.1}
\sup_{\mathbf 0 \neq \mathbf u_h \in \mathbf V^h_0(\Omega_{\mathrm
S})} \frac{(\mathrm{div}\,\mathbf u_h, q_h)_{\Omega_{\mathrm
S}}}{\|\mathbf u_h\|_{\mathbf H^1(\Omega_{\mathrm S})}}\ge
\beta\| q_h\|_{\Omega_{\mathrm S}}\qquad \forall q_h \in
L^h_0(\Omega_{\mathrm S}),
\end{equation}
\begin{equation}\label{eq:2.2}
\sup_{\mathbf 0 \neq \mathbf u_h \in \mathbf H^h_0(\Omega_{\mathrm
D})} \frac{(\mathrm{div}\,\mathbf u_h, q_h)_{\Omega_{\mathrm
D}}}{\|\mathbf u_h\|_{\mathbf H(\mathrm{div},\Omega_{\mathrm
D})}}\ge \beta\| q_h\|_{\Omega_{\mathrm D}}\qquad
\forall q_h \in L^h_0(\Omega_{\mathrm D}),
\end{equation}
\begin{equation}\label{eq:2.2b}
\mathrm{div}\,\mathbf
H^h(\Omega_{\mathrm D}) \subset L^h(\Omega_{\mathrm D}),
\end{equation}
\begin{equation}\label{eq:2.3}
\exists \mathbf v_h \in \mathbf V^h_{\mathrm
S}(\Omega_{\mathrm S}) \quad \mbox{s.t.}\quad  
\langle\mathbf v_h\cdot\mathbf{n},1\rangle_\Sigma\ge \beta \quad \text{and}\quad 
 \|\mathbf v_h\|_{\mathbf H^1(\Omega_{\mathrm S})}\le \beta,
\end{equation}
\begin{equation}\label{eq:2.3b}
(\mathbf L_h\phi_h)\cdot\mathbf{n} \qquad 
\|\mathbf L_h\phi_h\|_{\mathbf H(\mathrm{div},\Omega_{\mathrm D})}\le \beta \|\phi\|_{H^{-1/2}(\Sigma)}
\qquad\forall \phi_h\in \Phi^h_{\mathrm D}.
\end{equation}
Then the discrete equations\eqref{eq:2.0} are uniquely solvable and the following error estimate holds:
\begin{multline*}
\|\mathbf u-\mathbf u_h\|_{\mathbb
X}+\|p-p_h\|_{\Omega_{\mathrm
D}}
\lesssim    \inf_{\mathbf
u^h_{\mathrm S} \in \mathbf V^h_{\mathrm S}(\Omega_{\mathrm S})} \|
\mathbf u_{\mathrm S}-\mathbf u^h_{\mathrm S}\|_{\mathbf
H^1(\Omega_{\mathrm S})} + \inf_{\mathbf u^h_{\mathrm D}\in \mathbf
H^h_{\mathrm D}(\Omega_{\mathrm D})} \| \mathbf u_{\mathrm
D}-\mathbf u^h_{\mathrm D}\|_{\mathbf
H(\mathrm{div},\Omega_{\mathrm D})}+ \\
\inf_{q_h \in \mathbb{Q}}\| p-q_h\|_\Omega  + \lambda(h) \, \left( \|p_{\mathrm{D}} - R^h_D p_{\mathrm{D}} \|_\Sigma + 
\|\mathbf u_{\mathrm S}\cdot\mathbf{n}-R^h_D(\mathbf
u_{\mathrm S}\cdot\mathbf{n})\|_{\Sigma}
\right).
\end{multline*}
Here $\lambda(h)\equiv 0$ if $\Phi_{\mathrm S}^h \subset \Phi_{\mathrm D}^h$ and 
$
\mathrm \lambda(h) \lesssim h^{1/2} 
$
otherwise.
\end{theorem}

Let us briefly discuss the five hypotheses in Theorem \ref{thm:main}. The inf-sup condition \eqref{eq:2.1} 
is necessary and sufficient for stability of the discretization of the Stokes equation with homogeneous 
boundary conditions. The inf-sup condition \eqref{eq:2.2} and the restriction \eqref{eq:2.2b} are standard 
conditions for stability of the discretization of the Darcy equations with homogeneous boundary condition 
on the normal trace. 

Condition \eqref{eq:2.3b} is the existence of a uniformly bounded right-inverse of the operator 
$\mathbf H^h_{\mathrm D}(\Omega_{\mathrm D})\ni \mathbf v_h \mapsto \mathbf v_h\cdot\mathbf{n} 
\in \Phi_h^{\mathrm D}$. As discussed in \cite[Section 5]{anteriorPaper}, this condition is satisfied 
for  Brezzi-Douglas-Marini  (BDM) and Raviart-Thomas (RT) elements (see below for their definitions) 
on general shape-regular triangulations in the plane, and on tetrahedrizations of the space that are 
quasi-uniform near the boundary $\Sigma$. Existence of $\mathbf L_h$ satisfying \eqref{eq:2.3b} for 
BDM and RT elements in general tetrahedrizations is an open question. Hypothesis \eqref{eq:2.3} is 
a very mild condition demanding that the discrete space for the Stokes condition can provide 
non-trivial flow in to the Darcy domain without a blow-up of the velocity field. This 
condition is discussed in \cite[Section 6]{anteriorPaper}, where it is shown that as long 
as the Stokes velocity space contains piecewise linear 
functions, this condition is satisfied.

\paragraph{Some examples}
For precise descriptions of the finite element spaces below, the reader is referred to \cite{BreFor91}, 
\cite{ErnGue04} and \cite{GirRav86}.
All choices below will be given with the following assumptions:
\begin{itemize}
\item Hypothesis \eqref{eq:2.3} will be assumed to hold.
\item The convergence orders of the Stokes and Darcy elements are chosen to match.
\item If the discretization is conforming, we will assume that  the Darcy partition $\Sigma^h_{\mathrm D}$ is either 
equal to or a refinement of $\Sigma^h_{\mathrm S}$.
\end{itemize}
The Brezzi-Douglas-Marini (sometimes called Brezzi-Douglas-Dur\'an-Fortin in the three dimensional case) 
is the mixed finite element that uses the spaces
\begin{eqnarray*}
\mathbf H^h(\Omega_{\mathrm D})&:=&\{\mathbf u_h \in \mathbf H(\mathrm{div},\Omega_{\mathrm D})\,:\,\mathbf u_h|_T\in 
\mathbb P_{k}(T)^d \quad\forall T\in \mathcal T^h_{\mathrm D}\},\\
 L^h(\Omega_{\mathrm D})&:=&\{ p_h :\Omega_{\mathrm D}\to \mathbb R\,:\, p_h|_T \in \mathbb P_{k-1}(T) \quad\forall 
T\in \mathcal T^h_{\mathrm D}\},
\end{eqnarray*}
for $k\ge 1$. We will refer to it as the BDM($k$) element. The BDM(1) element can be coupled in a conforming way 
with the MINI element and the Bernardi-Raugel element. It can also be coupled with the $\mathbb P_2$-iso-$\mathrm P_1$ 
element in a conforming way if $\Sigma^h_{\mathrm D}$ is either equal to or a refinement of $\Sigma^{h/2}_{\mathrm S}$ 
and in a non-conforming way otherwise. The BDM(2) element can be coupled in a conforming way with the conforming 
Crouzeix-Raviart element. More generally speaking, BDM($k$) can be coupled with the Taylor-Hood element of order 
$k$ for any $k\ge 2$.

The Raviart-Thomas element of order $k$, henceforth referred to as RT($k$), is defined as the pair
\begin{eqnarray*}
\mathbf H^h(\Omega_{\mathrm D})&:=&\{\mathbf u_h \in \mathbf H(\mathrm{div},\Omega_{\mathrm D})\,:\,\mathbf u_h|_T\in 
\mathrm{RT}_k(T) \quad\forall T\in \mathcal T^h_{\mathrm D}\},\\
 L^h(\Omega_{\mathrm D})&:=&\{ p_h :\Omega_{\mathrm D}\to \mathbb R\,:\, p_h|_T \in \mathbb P_{k}(T) \quad\forall T\in 
\mathcal T^h_{\mathrm D}\},
\end{eqnarray*}
where $\mathrm{RT}_k(T)=\{\mathbf p(\mathbf x)+q(\mathbf x)\,\mathbf x\,:\, \mathbf p\in \mathbb P_k(T)^d, 
\quad q\in \mathbb P_k(T)\}$. The RT(0) element can be coupled in non-conforming way with the MINI element 
and the Bernardi-Raugel element. For $k\ge1 $, RT($k-1$) can be coupled with the Taylor-Hood element of order $k$.

\section{An alternative point of view}\label{sec:5}

In this section we propose a different way of interpreting the  coupled method, based on seeing 
the Darcy equations as part of a generalized boundary condition for the Stokes problem.

\subsection{The Darcy boundary condition}

Given $f_{\mathrm D}$ satisfying the compatibility condition \eqref{eq:0}, we consider the solution of the Darcy problem
\begin{eqnarray*}
\mathbf K^{-1}\mathbf u_{\mathrm D}^f+\nabla p_{\mathrm D}^f=\mathbf 0 & & \mbox{in $\Omega_{\mathrm D}$},\\
\mathrm{div}\,\mathbf u_{\mathrm D}^f = f_{\mathrm D}& & \mbox{in $\Omega_{\mathrm D}$},\\
\mathbf u_{\mathrm D}^f\cdot\mathbf{n} = 0 & & \mbox{on $\Gamma_{\mathrm D}\cup \Sigma$},\\
\int_{\Omega_{\mathrm D}} p_{\mathrm D}^f=0,
\end{eqnarray*}
and note that $p_{\mathrm D}^f\in H^1(\Omega_{\mathrm D})$. Also, for $\phi\in L^2(\Sigma)$, we consider the solution of
\begin{eqnarray*}
\mathbf K^{-1}\mathbf u_{\mathrm D}^\phi+\nabla p_{\mathrm D}^\phi=\mathbf 0 & & \mbox{in $\Omega_{\mathrm D}$},\\
\mathrm{div}\,\mathbf u_{\mathrm D}^\phi =\frac1{|\Omega|}\int_\Sigma \phi& & \mbox{in $\Omega_{\mathrm D}$},\\
\mathbf u_{\mathrm D}^\phi\cdot\mathbf{n} = 0 & & \mbox{on $\Gamma_{\mathrm D}$},\\
\mathbf u_{\mathrm D}^\phi\cdot\mathbf{n} = \phi & & \mbox{on $\Sigma$},\\
\int_{\Omega_{\mathrm D}} p_{\mathrm D}^\phi=0,
\end{eqnarray*}
and define with it the Flux-to-Pressure operator $\mathrm{FtP}(\phi):=p_{\mathrm D}^\phi|_\Sigma$. It is 
simple to prove that $\mathrm{FtP}$ is a linear and symmetric operator in $L^2(\Sigma)$. Indeed, the Flux-to-Pressure 
operator satisfies
\begin{equation}\label{eq:7.1}
\langle\mathrm{FtP}(\phi),\mathbf v\cdot\mathbf{n}\rangle_\Sigma=a_{\mathrm{D}}(
\mathbf u_{\mathrm D}^\phi,\mathbf v) - b_{\mathrm{D}}(\mathbf v,p_{\mathrm D}^\phi)
\end{equation}
for all $\mathbf v\in \mathbf H_{\mathrm D}(\mathrm{div},\Omega_{\mathrm D})$ such that 
$\mathbf{v}\cdot\mathbf{n}\in L^2(\Sigma)$, which gives
\begin{eqnarray}\label{sym}
\langle \mathrm{FtP}(\phi),\psi\rangle_\Sigma=\langle \mathrm{FtP}(\phi),\mathbf 
u_{\mathrm D}^\psi\cdot\mathbf{n}\rangle_\Sigma=a_{\mathrm{D}}(
\mathbf u_{\mathrm D}^\phi,\mathbf 
u_{\mathrm D}^\psi)
= \langle \mathrm{FtP}(\psi),\phi\rangle_\Sigma \quad \forall \phi,\psi\in L^2(\Sigma).\label{eq:7.2}
\end{eqnarray}

It is clear that, as $\int_\Sigma \mathbf u_{\mathrm S}\cdot\mathbf{n}=0$, we have the splitting 
\[
 p_{\mathrm D}|_\Sigma = p_{\mathrm D}^f+\mathrm{FtP}(\mathbf u_{\mathrm S}\cdot\mathbf{n})
\]
for the Darcy pressure on $\Sigma$. This allow us to write the coupling conditions 
\eqref{eq:00.4} as a unilateral 
boundary condition for the Stokes flow on the interface $\Sigma$:
\begin{equation}\label{eq:7.3}
2\nu\boldsymbol\varepsilon(\mathbf u_{\mathrm S})\mathbf{n}-p_{\mathrm S}\mathbf{n}+ 
\underbrace{\kappa\boldsymbol\pi_t\mathbf u_{\mathrm S}+\mathrm{FtP}(\mathbf u_{\mathrm S}
\cdot\mathbf{n})\mathbf{n}}=-p_{\mathrm D}^f \mathbf{n}.
\end{equation}
The underbracketed term corresponds to a symmetric positive semidefinite non-local operator that takes 
into account the influence of the Darcy domain on the Stokes flow, acting  separately on the tangential 
and normal components of the Stokes flow. The Stokes system \eqref{eq:00.1} 
can then be complemented with the non-local condition \eqref{eq:7.3} and the Dirichlet condition on 
$\Gamma_{\mathrm S}$ (see \eqref{eq:00.3})  to produce a formulation of the 
Stokes-Darcy problem that is equivalent to \eqref{eq:1.vp}. It consists in  
looking for $\mathbf u_{\mathrm S}\in \mathbf H^1_{\mathrm S}(\Omega_{\mathrm S})$ and $p_{\mathrm S}
\in L^2(\Omega_{\mathrm S})$ such that
\begin{equation}\label{eq:7.4}
\begin{array}{rll}
a_{\mathrm{S}}(\mathbf u_{\mathrm S},\mathbf v) + 
c(\mathbf u_{\mathrm S},\mathbf v) - b_{\mathrm{S}}(\mathbf v,p_{\mathrm S}) &= 
\ell(\mathbf v)\quad & \forall \mathbf v \in \mathbf H^1_{\mathrm S}(\Omega_{\mathrm S}),\\[2ex]
b_{\mathrm{S}}(\mathbf u_{\mathrm S},q)&=0 \qquad & \forall q\in L^2(\Omega_{\mathrm S}),
\end{array}
\end{equation}
where
\begin{equation}\label{eq:7.6}
c(\mathbf u,\mathbf v):= \langle\mathrm{FtP}(\mathbf u\cdot\mathbf{n}),\mathbf v\cdot\mathbf{n}\rangle_\Sigma
\end{equation}
and
\[
\ell (\mathbf v):= (\mathbf f_{\mathrm S},\mathbf v)_{\Omega_{\mathrm S}}-\langle p_{\mathrm D}^f,
\mathbf v\cdot\mathbf{n}\rangle_\Sigma.
\]
By \eqref{eq:7.2}, it follows that the bilinear form in \eqref{eq:7.6} is symmetric and positive 
semidefinite.  A simple argument shows that the bilinear form 
$a_{\mathrm{S}}(\mathbf u_{\mathrm S},\mathbf v) + 
c(\mathbf u_{\mathrm S},\mathbf v)$ is coercive in $\mathbf H^1(\Omega_{\mathrm S})$. This fact gives a very 
simple proof of the fact that the Stokes-Darcy system is well posed and that it can be understood 
as a modified Stokes problem without losing any of its good properties. This will be exploited to 
design an effective Krylov-based iterative method to solve the algebraic linear system 
of equations arising from the discrete counterpart of \eqref{eq:7.4}.

\subsection{The discrete flux-to-pressure operator}

If we now choose discrete spaces for the Darcy problem satisfying \eqref{eq:2.2}-\eqref{eq:2.2b}, 
we can define a discrete version of the operator $\mathrm{FtP}$ as follows. Given $\phi_h\in 
\Phi^h_{\mathrm D}$ with $\int_\Sigma \phi_h = 0$, we define $\mathrm{FtP}_h(\phi_h):\Phi^h_{\mathrm D}\to \mathbb R$ to be 
the functional (compare with \eqref{eq:7.1})
\begin{equation}\label{action}
\langle \mathrm{FtP}_h(\phi_h),\mathbf v_h\cdot\mathbf{n}\rangle_\Sigma:=a_{\mathrm D}(
\mathbf u_h^\phi,\mathbf v_h) - b_{\mathrm{D}}(\mathbf v_h,p_h^\phi),\quad \forall 
\mathbf v_h\in \mathbf H^h_{\mathrm D}(\Omega_{\mathrm D}),
\end{equation}
where $(\mathbf u_h^\phi,p_h^\phi)\in \mathbf H^h_{\mathrm D}(\Omega_{\mathrm D})\times 
L^h_0(\Omega_{\mathrm D})$ solves the discrete equations:
\begin{eqnarray}\label{darcy}
\mathbf u_h\cdot \mathbf{n} = \phi_h & & \mbox{on $\Sigma$},\\[1ex]
\mathbf u_h\cdot \mathbf{n} = 0 & & \mbox{on $\Gamma_{\mathrm D}$},\\[1ex]
a_{\mathrm D}(\mathbf u_h^\phi,\mathbf v_h) - b_{\mathrm{D}}(
\mathbf v_h,p_h^\phi)=0 & &  \forall \mathbf v_h \in \mathbf H^h_0(\Omega_{\mathrm D}),\\[1ex]
b_{\mathrm{D}}(\mathbf u_h^\phi,q_h)=0 & & \forall q_h\in L^h_0(\Omega_{\mathrm D}).
\end{eqnarray}
With arguments similar to those used in the continuous case, it is easy to prove that
\[
\langle \mathrm{FtP}_h(\phi_h),\psi_h\rangle_\Sigma = a_{\mathrm D}(\mathbf u_h^\phi,\mathbf 
u_h^\psi)=\langle \mathrm{FtP}_h(\psi_h),\phi_h\rangle_\Sigma\qquad 
\forall\phi_h,\psi_h \in \Phi^h_{\mathrm D},
\]
which shows that  the discrete flux-to-pressure operator $\mathrm{FtP}_h$
is also symmetric and nonnegative. 

The discrete pressure due to sources, $\gamma_h^f$, can be similarly defined as a residual:
\[
\langle \gamma_h^f,\mathbf v_h\cdot\mathbf{n}\rangle_\Sigma:=a_{\mathrm D}(\mathbf u_h^f,
\mathbf v_h) - b_{\mathrm D}(\mathbf v_h,p_h^f),
\]
where $(\mathbf u_h^f,p_h^f)\in \mathbf H^h_0(\Omega_{\mathrm D})\times L^h_0(\Omega_{\mathrm D})$ 
solve the discrete equations:
\begin{eqnarray*}
a_{\mathrm D}(\mathbf u_h^f,\mathbf v_h) - b_{\mathrm D}(\mathbf v_h,p_h^f)&=&0  \qquad \qquad\forall 
\mathbf v_h \in \mathbf H^h_0(\Omega_{\mathrm D}),\\
b_{\mathrm D}(\mathbf u_h^f,q_h)&=&(f_{\mathrm D},q_h)_{\Omega_{\mathrm D}} 
\qquad\forall q_h\in L^h_0(\Omega_{\mathrm D}).
\end{eqnarray*}

We recall that the operator $R^h$ is the $L^2(\Sigma)$-projection 
on $\Phi^h_{\mathrm D}$. It is straightforward that the discrete Darcy pressure and velocity of  
problem \eqref{eq:2.0} admit the splitting 
\[
 p_{\mathrm{D}}^h = p_h^f + p_h^{R^h_D(\mathbf u^h_{\mathrm S}\cdot\mathbf{n})}\quad \text{and} \quad 
 \mathbf{u}_{\mathrm{D}}^h = \mathbf{u}_h^f + \mathbf{u}_h^{R^h_D(\mathbf u^h_{\mathrm S}\cdot\mathbf{n})}.
\]
It follows that \eqref{eq:2.0} may be equivalently stated as follows: find $(\mathbf u^h_{\mathrm S},
p^h_{\mathrm S})\in \mathbf V^h_{\mathrm S}(\Omega_{\mathrm S})\times L^h(\Omega_{\mathrm S})$ such that
\begin{equation}\label{eq:7.7}
\begin{array}{rll}
a_{\mathrm{S}}(\mathbf u^h_{\mathrm S},\mathbf v_h)
+ c_h(\mathbf u^h_{\mathrm S},\mathbf v_h) - b_{\mathrm{S}}(\mathbf v_h,p^h_{\mathrm S}) 
&= \ell_h(\mathbf v_h)\quad &\forall \mathbf v_h \in \mathbf V^h_{\mathrm S}(\Omega_{\mathrm S}),\\[2ex]
b_{\mathrm{S}}(\mathbf u^h_{\mathrm S},q_h)&=0 \qquad &\forall q_h\in L^h(\Omega_{\mathrm S}),
\end{array}
\end{equation}
where
\[
c_h(\mathbf u_h,\mathbf v_h):= \langle\mathrm{FtP}_h(R^h(\mathbf u_h
\cdot\mathbf{n})),R^h(\mathbf v_h\cdot\mathbf{n})\rangle_\Sigma
\]
and
\[
\ell_h(\mathbf v_h):= (\mathbf f_{\mathrm S},\mathbf v_h)_{\Omega_{\mathrm S}}-
\langle p_{h}^f,R^h(\mathbf v_h\cdot\mathbf{n})\rangle_\Sigma.
\]
Inn the conforming case ($\Phi^h_{\mathrm S}\subset \Phi^h_{\mathrm D}$), the $L^2(\Sigma)$-projection 
operator $R^h$ does not play any role in the formulation. 

\section{The decoupled iterative method}\label{sec:6}
We introduce the self-adjoint operators  
$A^h_{\mathrm{S}}$ and $C^h_{\mathrm{S}}$ defined from $\mathbf V^h_{\mathrm S}(\Omega_{\mathrm S})$ to its 
dual $\mathbf V^h_{\mathrm S}(\Omega_{\mathrm S})^*$ by 
\[
 \langle A^h_{\mathrm{S}}\mathbf{u} ,  \mathbf{v}\rangle_{\mathbf V^h_{\mathrm S}(\Omega_{\mathrm S})^* 
 \times \mathbf V^h_{\mathrm S}(\Omega_{\mathrm S})} = a_{\mathrm{S}}(\mathbf{u}, \mathbf{v}) 
 \quad \text{and} \quad 
 \langle C^h_{\mathrm{S}}\mathbf{u} ,  \mathbf{v}\rangle_{\mathbf V^h_{\mathrm S}(\Omega_{\mathrm S})^* 
 \times \mathbf V^h_{\mathrm S}(\Omega_{\mathrm S})} = c_h(\mathbf{u}, \mathbf{v})
\]
for all $\mathbf{u}, \mathbf{v}\in \mathbf V^h_{\mathrm S}(\Omega_{\mathrm S})$. Let also
$B^h_{\mathrm{S}}:\, \mathbf V^h_{\mathrm S}(\Omega_{\mathrm S}) \to L^h(\Omega_{\mathrm S})^*$ be the operator 
defined by 
\[
 \langle B^h_{\mathrm{S}} \mathbf{u}, q 
 \rangle_{L^h(\Omega_{\mathrm S})^*\times  L^h(\Omega_{\mathrm S})} = -b_{\mathrm{S}}(\mathbf u,q)
\]
for all $\mathbf{u}\in \mathbf V^h_{\mathrm S}(\Omega_{\mathrm S})$ and $q\in L^h(\Omega_{\mathrm S})$.  

Problem \eqref{eq:7.7} can be written in operator form as follows:
\begin{equation}\label{matrix1}
 \mathcal{A}^h_{\mathrm{S}}
\begin{pmatrix}
 \mathbf{u}^h_{\mathrm{S}}\\p^h_{\mathrm{S}}
\end{pmatrix}
=
\begin{pmatrix}
 \ell_h\\ 0
\end{pmatrix}
\end{equation}
where
\[
 \mathcal{A}^h_{\mathrm{S}}:=
\begin{pmatrix}
 A^h_{\mathrm{S}}+C^h_{\mathrm{S}} & (B^h_{\mathrm{S}})^{\mathtt t}\\
 B^h_{\mathrm{S}} & \mathbf{0}
\end{pmatrix}:\, \mathbf{V}^h_{\mathrm S}(\Omega_{\mathrm S})\times  L^h(\Omega_{\mathrm S}) \to 
\mathbf{V}^h_{\mathrm S}(\Omega_{\mathrm S})^*\times  L^h(\Omega_{\mathrm S})^*
\]
and $(B^h_{\mathrm{S}})^{\mathtt t}$ is the adjoint of $B^h_{\mathrm{S}}$. We know from Theorem 
\ref{thm:main} that both $\|\mathcal{A}^h_{\mathrm{S}}\|$ and $\|(\mathcal{A}^h_{\mathrm{S}})^{-1}\|$ are uniformly bounded in 
$h$. If we denote by $I^h_{\mathrm S}:L^h(\Omega_{\mathrm S})\to L^h(\Omega_{\mathrm S})^*$ the Riesz operator defined by
\[
 \langle I^h_{\mathrm{S}}p ,  q\rangle_{ L^h(\Omega_{\mathrm S})^* 
 \times L^h(\Omega_{\mathrm S}) } = (p, q)_{\Omega_{\mathrm S}}
 \quad \forall p,q \in L^h(\Omega_{\mathrm S}),
\]
then, 
the positive-definite self-adjoint operator 
\[
 \mathcal{P}^h_{\mathrm{S}}:= \begin{pmatrix}
 A^h_{\mathrm{S}} & 0\\
 0 & I^h_{\mathrm S}
\end{pmatrix}:\, \mathbf{V}^h_{\mathrm S}(\Omega_{\mathrm S}) \times L^h(\Omega_{\mathrm S}) \to 
\mathbf{V}^h_{\mathrm S}(\Omega_{\mathrm S})^* \times L^h(\Omega_{\mathrm S})^*
\]
and its inverse are uniformly bounded uniformly in $h$. 
It follows that the condition number of 
$(\mathcal{P}^h_{\mathrm{S}})^{-1}\mathcal{A}^h_{\mathrm{S}}$ is bounded from above by a constant 
independent of the mesh parameter $h$.
Consequently, the MINRES algorithm preconditioned with $(\mathcal{P}^h_{\mathrm{S}})^{-1}$ solves  
\eqref{matrix1}  with a reduction of the norm of the residual that is  
independent of the mesh size $h$.

Let us now discuss how  action of  $C^h_{\mathrm{S}}$ on a given 
$\mathbf{u}_{\mathrm{S}}^h\in \mathbf{V}^h_{\mathrm S}(\Omega_{\mathrm S})$. 
To this end we introduce the self-adjoint operators  
$A^h_{\mathrm{D}}$ and $D^h_{\mathrm{D}}$ defined from $\mathbf H^h_0(\Omega_{\mathrm D})$ to its 
dual $\mathbf H^h_0(\Omega_{\mathrm D})^*$ by 
\[
 \langle A^h_{\mathrm{D}}\mathbf{u} ,  \mathbf{v}\rangle_{\mathbf H^h_{\mathrm D}(\Omega_{\mathrm D})^* 
 \times \mathbf H^h_{\mathrm D}(\Omega_{\mathrm D})} = a_{\mathrm{D}}(\mathbf{u}, \mathbf{v}) 
 \quad \text{and} \quad 
 \langle D^h_{\mathrm{D}}\mathbf{u} ,  \mathbf{v}\rangle_{\mathbf H^h_{\mathrm D}(\Omega_{\mathrm D})^* 
 \times \mathbf H^h_{\mathrm D}(\Omega_{\mathrm D})} = (\div\mathbf{u}, \div \mathbf{v})_{\Omega_{\mathrm{D}}}
\]
for all $\mathbf{u}, \mathbf{v}\in \mathbf H^h_0(\Omega_{\mathrm D})$. Let us also consider 
$B^h_{\mathrm{D}}:\, \mathbf H^h_0(\Omega_{\mathrm D}) \to L_0^h(\Omega_{\mathrm D})^*$ given by 
\[
 \langle B^h_{\mathrm{D}} \mathbf{u}, q 
 \rangle_{L_0^h(\Omega_{\mathrm D})^*\times  L_0^h(\Omega_{\mathrm D})} = -b_{\mathrm{D}}(\mathbf u,q)
\]
for all $\mathbf{u}\in \mathbf H^h_0(\Omega_{\mathrm D})$ and $q\in L_0^h(\Omega_{\mathrm D})$.

We compute $C^h_{\mathrm{S}}\mathbf{u}_{\mathrm{S}}^h$ 
through \eqref{action} after solving problem \eqref{darcy} with 
$\phi_h = R_h(\mathbf{u}_{\mathrm{S}}^h\cdot \mathbf{n})$. This is to say that 
we have to deal with a saddle point problem of the form
\begin{equation}\label{matrix2}
 \mathcal{A}^h_{\mathrm{D}}
\begin{pmatrix}
 \mathbf{u}_h^\phi\\p^\phi_h
\end{pmatrix}
=
\begin{pmatrix}
 \mathcal{F}_h\\ \mathcal{G}_h
\end{pmatrix}
\end{equation}
where
\[
 \mathcal{A}^h_{\mathrm{D}}:=
\begin{pmatrix}
 A^h_{\mathrm{D}} & (B^h_{\mathrm{D}})^{\mathtt t}\\
 B^h_{\mathrm{D}} & \mathbf{0}
\end{pmatrix}:\, \mathbf{H}^h_0(\Omega_{\mathrm D})\times  L_0^h(\Omega_{\mathrm D}) \to 
\mathbf{H}^h_0(\Omega_{\mathrm D})^*\times  L_0^h(\Omega_{\mathrm D})^*
\]
and $(B^h_{\mathrm{D}})^{\mathtt t}$ is the adjoint of $B^h_{\mathrm{D}}$.

The stability of the pair of spaces 
$(\mathbf H^h_0(\Omega_{\mathrm D}), 
L^h_0(\Omega_{\mathrm S}))$ \eqref{eq:2.2}-\eqref{eq:2.2b} ensures that 
both $\|\mathcal{A}^h_{\mathrm{D}}\|$ and $\|(\mathcal{A}^h_{\mathrm{D}})^{-1}\|$ are uniformly bounded in 
$h$. If we denote by $I^h_{\mathrm D}:\: L_0^h(\Omega_{\mathrm D}) \to L_0^h(\Omega_{\mathrm D})^*$ 
the Riesz operator given by 
\[
 \langle I^h_{\mathrm{D}}p ,  q\rangle_{ L_0^h(\Omega_{\mathrm D})^* 
 \times L_0^h(\Omega_{\mathrm D}) } = (p, q)_{\Omega_{\mathrm D}}
 \quad \forall p,q \in L_0^h(\Omega_{\mathrm D}),
\]
it is clear that 
the block diagonal positive-definite self-adjoint operator 
\[
 \mathcal{P}^h_{\mathrm{D}}:= \begin{pmatrix}
 A^h_{\mathrm{D}} + D^h_{\mathrm{D}} & 0\\
 0 & I^h_{\mathrm D}
\end{pmatrix}:\, \mathbf{H}^h_0(\Omega_{\mathrm D}) \times L_0^h(\Omega_{\mathrm D}) \to 
\mathbf{H}^h_0(\Omega_{\mathrm D})^* \times L_0^h(\Omega_{\mathrm D})^*
\]
and its inverse are also uniformly bounded in $h$. 
It follows that we can find an inclusion set for the eigenvalues of 
$(\mathcal{P}^h_{\mathrm{D}})^{-1}\mathcal{A}^h_{\mathrm{D}}$ that is  
independent of $h$. This means that the MINRES method preconditioned with $(\mathcal{P}^h_{\mathrm{D}})^{-1}$
yields the solution of \eqref{matrix2}  in a number of iterations independent on the 
mesh size $h$.

Summing up, the decoupled iterative method we are proposing here 
to solve \eqref{eq:2.0} consists in  two nested MINRES 
algorithms. 
Computationally, the actions of the preconditioners correspond to solving two decoupled local problems. 
The first one is defined by the bilinear form $a_{\mathrm{S}}(\cdot, \cdot)$ in $\mathbf{V}^h_{\mathrm{S}}(\Omega_{\mathrm{S}})$  
and corresponds to the block $A^h_{\mathrm{S}}$.  Actually, $A^h_{\mathrm{S}}$ is associated with the operator 
$-2\nu\textbf{div}(\boldsymbol\varepsilon(\cdot))$. Therefore, the local problem 
in the fluid amounts to a vector Laplace equation with a Dirichlet boundary condition on $\Gamma_{\mathrm{S}}$,  
a Neumann condition in the normal direction and the slip boundary condition in the tangential 
direction on $\Sigma$. The other local problem is defined by the bilinear form 
\[
 (\mathbf{K}^{-1}\mathbf{u}_{\mathrm{D}}, \mathbf{v}_{\mathrm{D}})_{\Omega_{\mathrm{D}}} + 
 (\div\mathbf{u}_{\mathrm{D}}, \div\mathbf{v}_{\mathrm{D}})_{\Omega_{\mathrm{D}}}
\]
on $\mathbf{H}_0^h(\Omega_{\mathrm{D}})$ corresponding to the  
diagonal block $A^h_{\mathrm{D}}+D^h_{\mathrm{D}}$.  

For the construction of practical preconditioners for discrete systems, the computational cost
of evaluating these operators and the memory requirements of these procedures are key factors.
The exact inverses appearing in the canonical preconditioners should be replaced by proper cost 
effective, and norm equivalent operators. Let us consider self-adjoint and positive-definite operators 
$P_\mathrm{S}^h$ and 
$P_\mathrm{D}^h$ that are spectrally equivalent to $A_{\mathrm{S}}^h$ and 
$A^h_{\mathrm{D}}+D^h_{\mathrm{D}}$ respectively,  i.e., 
\[
 \langle A_{\mathrm{S}}^h\mathbf{u}_{\mathrm{S}},\,  \mathbf{u}_{\mathrm{S}} 
 \rangle_{\mathbf{H}_{\mathrm{S}}(\Omega_{\mathrm{S}})^*\times \mathbf{H}_{\mathrm{S}}(\Omega_{\mathrm{S}})} \simeq 
 \langle P_{\mathrm{S}}^h\mathbf{u}_{\mathrm{S}},\,  \mathbf{u}_{\mathrm{S}} 
 \rangle_{\mathbf{H}_{\mathrm{S}}(\Omega_{\mathrm{S}})^*\times \mathbf{H}_{\mathrm{S}}(\Omega_{\mathrm{S}})} 
\] 
and
\[
 \langle (A^h_{\mathrm{D}}+D^h_{\mathrm{D}})\mathbf{u}_{\mathrm{D}},\,  \mathbf{u}_{\mathrm{D}} 
 \rangle_{\mathbf{H}_0(\Omega_{\mathrm{D}})^*\times \mathbf{H}_0(\Omega_{\mathrm{D}})} \simeq
 \langle P_{\mathrm{D}}^h\mathbf{u}_{\mathrm{D}},\,  \mathbf{u}_{\mathrm{D}} 
 \rangle_{\mathbf{H}_0(\Omega_{\mathrm{D}})^*\times \mathbf{H}_0(\Omega_{\mathrm{D}})}.
\]
Then, instead of 
$(\mathcal{P}^h_{\mathrm{S}})^{-1}$ and $(\mathcal{P}^h_{\mathrm{D}})^{-1}$, 
we can  use respectively the preconditioners
\[
 \begin{pmatrix}
 (P^h_{\mathrm{S}})^{-1} & 0\\
 0 & (I^h_{\mathrm S})^{-1}
\end{pmatrix}
\quad \text{and} \quad 
 \begin{pmatrix}
 (P^h_{\mathrm{D}})^{-1} & 0\\
 0 & (I^h_{\mathrm D})^{-1}
\end{pmatrix}
\]
and still have an optimal decoupled iterative method for problem \eqref{eq:2.0}. 
Ideally, we would have the actions 
of $(P^h_{\mathrm{S}})^{-1}$ and $(P^h_{\mathrm{D}})^{-1}$ cost about the same as the actions of 
$A_{\mathrm{S}}^h$ and $A^h_{\mathrm{D}}+D^h_{\mathrm{D}}$. 
As $A_{\mathrm{S}}^h$ corresponds to a second-order elliptic 
operator in $\mathbf{H}_{\mathrm{S}}^1(\Omega_{\mathrm{S}})$, we can easily take advantage of multigrid 
techniques or domain decomposition methods to find a good candidate for $(P^h_{\mathrm{S}})^{-1}$. 
The construction of a preconditioner $(P^h_{\mathrm{D}})^{-1}$ is less obvious. 

\subsection{Nodal auxiliary space preconditioning in $\mathbf{H}(\div)$}
In this section we describe the construction of the 
nodal auxiliary space preconditioning of Hiptmair and Xu \cite{HiptmairXu} 
for elliptic problems in $\mathbf{H}_0(\div,\Omega_{\mathrm{D}})$. This is our choice 
here for the matrix version $(\mathbf{P}_{\mathrm{D}}^h)^{-1}$ of the preconditioner 
$(P^h_{\mathrm{D}})^{-1}$ needed in the last section. 
To fix the ideas, we assume that  the tensor $\mathbf{K}$ is given by 
$\tau^{-1}\mathbf{I}$ where $\mathbf{I}$ is  the identity in $\mathbb{R}^d$ and 
$\tau$ is a given positive constant. In our numerical experiments,  $\mathbf{H}^h_0(\Omega_{\mathrm{D}})$ 
is derived from the RT($k-1$) or BDM($k$) mixed finite elements with $k=1$  or 2. Let 
$\{\boldsymbol{\phi}_i; \quad i=1,\dots, I \}$ be the usual basis of $\mathbf{H}_0^h(\Omega_{\mathrm{D}})$, 
then the matrix realizations of $A^h_{\mathrm{D}}$ and $D^h_{\mathrm{D}}$ are given by
\[
 \mathbf{A}^h_{\mathrm{D}} :=  ( \tau (\boldsymbol{\phi}_i, \boldsymbol{\phi}_i )_{\Omega_{\mathrm{D}}} )
 _{1 \leq i, j\leq I}
\]
and 
\[
 \mathbf{D}^h_{\mathrm{D}} := ( (\div \boldsymbol{\phi}_i, \div \boldsymbol{\phi}_i)
 _{\Omega_{\mathrm{D}}})_{1 \leq i, j\leq I}
\]
respectively. Our aim is to 
provide a matrix $\mathbf{P}_{\mathrm{D}}^h\in \mathbb{R}^{I\times I}$  
that is spectrally equivalent to $\mathbf{A}^h_{\mathrm{D}} + \mathbf{D}^h_{\mathrm{D}}$
and such that the action $(\mathbf{P}_{\mathrm{D}}^h)^{-1}$ on a given vector is easier to compute then 
that of $(\mathbf{A}^h_{\mathrm{D}} + \mathbf{D}^h_{\mathrm{D}})^{-1}$. 

We denote by $[V^h_0(\Omega_{\mathrm{D}})]^d\subset \mathbf{H}_0^1(\Omega_{\mathrm{D}})$ the 
standard space of piecewise $\mathbb{P}_k$ and continuous vector fields and consider  its usual 
nodal basis $\{\boldsymbol{\varphi}_\ell,\, \ell=1,\cdots,  Ld\}$, where $L$ is the dimension of 
$V^h_0(\Omega_{\mathrm{D}})$. We introduce 
the matrix $\mathbf{L}_h$ given by
\[
 (\mathbf{L}_h)_{\ell,k} = (\nabla \boldsymbol{\varphi}_\ell, \nabla\boldsymbol{\varphi}_m)_{\Omega_{\mathrm{D}}} + \tau 
 (\boldsymbol{\varphi}_\ell, \boldsymbol{\varphi}_m)_{\Omega_{\mathrm{D}}}, \quad 1\leq \ell,m\leq Ld.
\]
 
In the three dimensional case ($d=3$), 
we also need to consider the N\'ed\'elec space $\mathbf{W}^h_0(\Omega_{\mathrm{D}}) \subset 
\mathbf{H}_0(\mathbf{curl}, \Omega_{\mathrm{D}})$ of order $k$. We denote its usual basis   
$\{\boldsymbol{\sigma}_i, \, i=1,\dots, N \}$. We introduce the diagonal matrix $\mathbf{S}_h^{\text{curl}}$ 
given by 
\[
(\mathbf{S}_h^{\text{curl}})_{i,i} := ( \textbf{curl}\, \boldsymbol{\sigma}_i, 
\textbf{curl}\, \boldsymbol{\sigma}_i )_{\Omega_{\mathrm{D}}}\quad i= 1,\cdots,N
\]
and denote the diagonal of $\mathbf{A}^h_{\mathrm{D}} + \mathbf{D}^h_{\mathrm{D}}$ by $\mathbf{S}_h^{\text{div}}$.

In the three-dimensional case, we denote by $\mathbf{C}_h\in \mathbb{R}^{N\times I}$ the matrix that represents 
$\mathbf{curl}:\, \mathbf{W}_0^h(\Omega_{\mathrm{D}})\to \mathbf{H}_0^h(\Omega_{\mathrm{D}})$ in the following sense,  
\[
 \mathbf{curl} \boldsymbol{\sigma}_i = \sum_{j= 1}^I (\mathbf{C}_h)_{i, j}\boldsymbol{\phi}_j, \quad \forall 
 i=1,\cdots, N.  
\]
In the the two-dimensional case, the matrix $\mathbf{C}_h\in \mathbb{R}^{L \times I}$ is defined similarly 
with respect to the operator $\mathbf{curl}:\, V^h_0(\Omega_{\mathrm{D}}) \to \mathbf{H}_0^h(\Omega_{\mathrm{D}})$ 
defined by $\mathbf{curl}v = \begin{pmatrix}\partial_2 v \\ -\partial_1 v \end{pmatrix}$.

We use  $\Pi_h^{\text{curl}}$ and $\Pi_h^{\text{div}}$ to denote the canonical 
interpolation operators onto the finite element spaces  
$\mathbf{W}^h_0(\Omega_{\mathrm{D}})$ and $\mathbf{H}_0^h(\Omega_{\mathrm{D}})$ respectively.
The mappings $\Pi_h^{\text{div}}:\,[V^h_0(\Omega_{\mathrm{D}})]^d \to  \mathbf{H}_0^h(\Omega_{\mathrm{D}})$  
and $\Pi_h^{\text{curl}}:\,[V^h_0(\Omega_{\mathrm{D}})]^3 \to  \mathbf{W}_0^h(\Omega_{\mathrm{D}})$ will be 
described by the matrices $\mathbf{I}_h^{\text{div}}\in \mathbb{R}^{dL \times I}$ ($d= 2,3$)
and $\mathbf{I}_h^{\text{curl}}\in \mathbb{R}^{3 L \times N}$ defined by 
\[
\Pi_h^{\text{div}}\,\boldsymbol{\varphi}_\ell  = \sum _{j= 1}^I (\mathbf{I}_h^{\text{div}})_{\ell, j} \boldsymbol{\phi}_j, 
\quad \forall \ell = 1,\cdots, Ld 
\]
and 
\[
\Pi_h^{\text{curl}}\, \boldsymbol{\varphi}_\ell  = \sum _{j= 1}^N (\mathbf{I}_h^{\text{curl}})_{\ell, j} \boldsymbol{\sigma}_j, 
\quad \forall \ell = 1,\cdots, 3 L,
\]
respectively.

The 3d-$\mathbf{H}(\div)$ auxiliary space preconditioner of Hiptmair and Xu consists in 
\[
 (\mathbf{P}_{\mathrm{D}}^h)^{-1} := (\mathbf{S}_h^{\text{div}})^{-1} +  
 \mathbf{I}_h^{\text{div}} (\mathbf{L}_h)^{-1} (\mathbf{I}_h^{\text{div}})^{\mathtt t} + 
 \tau^{-1}\mathbf{C}_h \Big( (\mathbf{S}_h^{\text{curl}})^{-1} 
 + 
 \mathbf{I}_h^{\text{curl}} (\mathbf{L}_h)^{-1} (\mathbf{I}_h^{\text{curl}})^{\mathtt t} 
 \Big) \mathbf{C}_h^{\mathtt t}
\]
and the 2-$d$ version of this preconditioner is defined by
\[
 (\mathbf{P}_{\mathrm{D}}^h)^{-1} := (\mathbf{S}_h^{\text{div}})^{-1} +  
 \mathbf{I}_h^{\text{div}} (\mathbf{L}_h)^{-1} (\mathbf{I}_h^{\text{div}})^{\mathtt t} + 
 \tau^{-1}\mathbf{C}_h \left( -\Delta_h \right)^{-1} \mathbf{C}_h^{\mathtt t}
\]
where the matrix $\Delta_h$ stands for the discrete Laplacian on the finite element space  
$V^h_0(\Omega_{\mathrm{D}})$.

Notice that the transfer matrices $\mathbf{C}_h$, $\mathbf{I}_h^{\text{div}}$  
and $\mathbf{I}_h^{\text{curl}}$ corresponding to the $\textbf{curl}$ operator  and the interpolations are sparse 
matrices that can be computed in a straightforward manner.  The evaluation the preconditioner 
is then essentially reduced to several second-order elliptic operators. Hence, standard multigrid techniques 
domain decomposition methods for $H^1$ equations can be applied.

\section{Numerical results}\label{sec:7}

This section is devoted to the description of numerical experiments validating the effectiveness of 
the decoupled iterative method. We will show results for two dimensional problems, considering three examples of pairs of stable elements for the Darcy-Stokes problem. 
The first two examples correspond to the conforming Galerkin schemes based on the combination of  
the MINI and $\mathbb P_2\mbox{-iso-}\mathbb P_1$  elements for the Stokes problem 
with the lowest order Brezzi-Douglas-Marini element BDM(1). 
The third one is the nonconforming scheme resulting from the Taylor-Hood element 
and the second order Raviart-Thomas element RT(1). 

\subsection{Convergence rates}
We begin by introducing some notation. 
The variable $DOF$ stands for the total number of degrees of freedom defining the 
finite element subspaces $\mathbb{X}^h$ and $\mathbb{Q}^h$, and 
the individual errors are denoted by:
\[
 \mathrm{e}(\mathbf{u}_{\mathrm D}) := \| \mathbf{u}_{\mathrm D} - \mathbf{u}^h_{\mathrm D} \|_{
\mathbf{H}(\mathrm{div}, \Omega_{\mathrm D})}, \qquad 
\mathrm{e}(\mathbf{u}_{\mathrm S}) := \| \mathbf{u}_{\mathrm S} - \mathbf{u}^h_{\mathrm S} \|_{
 \mathbf H^1(\Omega_{\mathrm S})
},
\]
and
 \[
  \mathrm{e}(p_{\mathrm D}) := \| p_{\mathrm D} - p^h_{\mathrm D}\|_{\Omega_{\mathrm D}},\qquad 
\mathrm{e}(p_{\mathrm S}) := \| p_{\mathrm S} - p^h_{\mathrm S}\|_{\Omega_{\mathrm S}},
 \]
where $\mathbf{u}^h_{\mathrm D} := \mathbf u_h|_{\Omega_{\mathrm D}}$, $\mathbf{u}^h_{\mathrm S} := \mathbf u_h|_{\Omega_{\mathrm D}}$, 
$p^h_{\mathrm D} := p_h|_{\Omega_{\mathrm D}}$ and $p^h_{\mathrm S} := p_h|_{\Omega_{\mathrm S}}$ with  
$(\mathbf u_h,p_h)\in \mathbb X^h\times \mathbb Q^h$ being the solution of \eqref{eq:2.0}. We also let 
$\verb"r"(\mathbf{u}_{\mathrm D})$, $\verb"r"(\mathbf{u}_{\mathrm S})$, $\verb"r"(p_{\mathrm D})$ and $\verb"r"(p_{\mathrm S})$ 
be the experimental rates of convergence given by 
\[
 \verb"r"(\mathbf{u}_{\mathrm D}) := \frac{\log(\mathrm{e}(\mathbf{u}_{\mathrm D})/ \mathrm{e}'(\mathbf{u}_{\mathrm D}) )}{\log(h/h')},\qquad
 \verb"r"(\mathbf{u}_{\mathrm S}) := \frac{\log(\mathrm{e}(\mathbf{u}_{\mathrm S})/ \mathrm{e}'(\mathbf{u}_{\mathrm S}) )}{\log(h/h')},
\]
and
\[
 \verb"r"(p_{\mathrm D}) := \frac{\log(\mathrm{e}(p_{\mathrm D})/ \mathrm{e}'(p_{\mathrm D}) )}{\log(h/h')},\qquad
 \verb"r"(p_{\mathrm S}) := \frac{\log(\mathrm{e}(p_{\mathrm S})/ \mathrm{e}'(p_{\mathrm S}) )}{\log(h/h')},
\]
where $h$ and $h'$ are two consecutive mesh sizes with errors $\mathrm{e}$ and $\mathrm{e}'$.

We now describe the data of the example. We consider the domains
$\Omega_\mathrm{D} := (0,1) \times (0,1/2)$ and $\Omega_\mathrm{S} := (0,1) \times (1/2,1)$.
We take $\nu =1$, $\kappa = 1$ and $\textbf{K} = \textbf{I}$, 
the identity of $\mathbb{R}^{2 \times 2}$.  The  right hand side 
functions are selected in the model in such a way 
that the exact solution is given by:
\[
p_\textrm{D}(\boldsymbol{x}):= 6\displaystyle\pi\left(\frac{x_2}{2}-\frac{1}{4\pi}\sin(2\pi x_2)\right)
\sin^2(2\pi x_1)\cos(2\pi x_1),
\]
in the porous media and by
$$
\textbf{u}_\textrm{S}(\boldsymbol{x}) := \displaystyle 2\pi \left(\begin{array}{c}
             \displaystyle \sin(\pi x_2)\cos(\pi x_2)\sin^3(2\pi x_1) \\[2ex]
             -3\sin^2(2\pi x_1)\cos(2\pi x_1)\sin^2(\pi x_2)
                                                 \end{array}
\right) 
$$
and
\[
p_\textrm{S}(\boldsymbol{x}) := -\displaystyle\frac{\pi}{4}\cos\left(\frac{\pi}{2}x_1\right)\left(
x_2+0.5-2\cos^2\left(\frac{\pi}{2}(x_2+0.5)\right)\right)
\]
in the Stokes domain.

\begin{table}[ht!]
\footnotesize
\begin{center}
\begin{tabular}{r| c ||c |c ||c |c ||c |c ||c |c}

$DOF$ & $h$ & $\mathrm{e}(\textbf{u}_\textrm{S})$ & $\verb"r"(\textbf{u}_\textrm{S})$ & $\mathrm{e}(p_\textrm{S})$ & $\verb"r"(p_\textrm{S})$ &
$\mathrm{e}(\textbf{u}_\textrm{D})$ & $\verb"r"(\textbf{u}_\textrm{D})$ & $\mathrm{e}(p_\textrm{D})$ & $\verb"r"(p_\textrm{D})$ \cr

\hline

543 & 1/8 & 1.86E$+$01 & $-$ & 9.26E$-$00 & $-$ & 4.73E$+$01 & $-$ & 1.60E$-$01 & $-$ \cr

2043 & 1/16 & 1.01E$+$01 & 0.87 & 3.04E$-$00 & 1.60 & 2.48E$+$01 & 0.92 & 8.10E$-$02 & 0.98 \cr

7923 & 1/32 & 5.17E$-$00 & 0.97 & 8.80E$-$01 & 1.79 & 1.26E$+$01 & 0.98 & 3.99E$-$02 & 1.02 \cr

31203 & 1/64 & 2.59E$-$00 & 0.99 & 2.49E$-$01 & 1.82 & 6.31E$-$00 & 0.99 & 1.98E$-$02 & 1.00 \cr

123843 & 1/128 & 1.29E$-$00 & 0.99 & 7.56E$-$02 & 1.72 & 3.16E$-$00 & 0.99 & 9.92E$-$03 & 1.00 \cr

\hline

\end{tabular}
\caption{Convergence rates: MINI--BDM(1)}
\label{table:3}
\end{center}
\end{table}

\begin{table}[ht!]
\footnotesize
\begin{center}
\begin{tabular}{r| c ||c |c ||c |c ||c |c ||c |c}

$DOF$ & $h$ & $\mathrm{e}(\textbf{u}_\textrm{S})$ & $\verb"r"(\textbf{u}_\textrm{S})$ & $\mathrm{e}(p_\textrm{S})$ & $\verb"r"(p_\textrm{S})$ &
$\mathrm{e}(\textbf{u}_\textrm{D})$ & $\verb"r"(\textbf{u}_\textrm{D})$ & $\mathrm{e}(p_\textrm{D})$ & $\verb"r"(p_\textrm{D})$ \cr

\hline

385 & 1/8 & 1.86E$+$01 & $-$ & 4.10E$-$00 & $-$ & 4.73E$+$01 & $-$ & 1.60E$-$01 & $-$ \cr

1423 & 1/16 & 1.01E$+$01 & 0.87 & 2.14E$-$00 & 0.94 & 2.48E$+$01 & 0.92 & 8.11E$-$02 & 0.98 \cr

5467 & 1/32 & 5.17E$-$00 & 0.97 & 6.03E$-$01 & 1.82 & 1.26E$+$01 & 0.98 & 3.99E$-$02 & 1.02 \cr

21427 & 1/64 & 2.59E$-$00 & 0.99 & 1.57E$-$01 & 1.94 & 6.32E$-$00 & 0.99 & 1.98E$-$02 & 1.00 \cr

84835 & 1/128 & 1.30E$-$00 & 0.99 & 4.25E$-$02 & 1.88 & 3.16E$-$00 & 0.99 & 9.92E$-$03 & 1.00 \cr

\hline

\end{tabular}
\caption{Convergence rates: $\mathbb{P}$2-iso-$\mathbb{P}$1--BDM(1)}
\label{table:4}
\end{center}
\end{table}

\begin{table}[ht!]
\footnotesize
\begin{center}
\begin{tabular}{r| c ||c |c ||c |c ||c |c ||c |c}

$DOF$ & $h$ & $\mathrm{e}(\textbf{u}_\textrm{S})$ & $\verb"r"(\textbf{u}_\textrm{S})$ & $\mathrm{e}(p_\textrm{S})$ & $\verb"r"(p_\textrm{S})$ &
$\mathrm{e}(\textbf{u}_\textrm{D})$ & $\verb"r"(\textbf{u}_\textrm{D})$ & $\mathrm{e}(p_\textrm{D})$ & $\verb"r"(p_\textrm{D})$ \cr

\hline

887 & 1/8 & 4.09E$-$00 & $-$ & 1.08E$-$00 & $-$ & 1.48E$+$01 & $-$ & 5.23E$-$02 & $-$ \cr

3371 & 1/16 & 9.56E$-$01 & 2.09 & 8.88E$-$02 & 3.60 & 4.03E$-$00 & 1.87 & 1.35E$-$02 & 1.95 \cr

13139 & 1/32 & 2.37E$-$01 & 2.01 & 7.06E$-$03 & 3.65 & 1.07E$-$00 & 1.90 & 3.40E$-$03 & 1.99 \cr

51875 & 1/64 & 5.93E$-$02 & 1.99 & 7.85E$-$04 & 3.17 & 2.94E$-$01 & 1.87 & 8.50E$-$04 & 2.00 \cr

206147 & 1/128 & 1.48E$-$02 & 1.99 & 1.98E$-$04 & 1.98 & 8.43E$-$02 & 1.80 & 2.12E$-$04 & 2.00 \cr

\hline

\end{tabular}
\caption{Convergence rates: Taylor-Hood--RT(1)}
\label{table:5}
\end{center}
\end{table}

We begin by providing a numerical exploration of the asymptotic 
convergence rates of the three examples.  
In Tables \ref{table:3}, \ref{table:4} and \ref{table:5}, we summarize the convergence 
history of the Galerkin scheme \eqref{eq:2.0} for a sequence of nested uniform meshes of the
computational domain $\Omega := (0,1)^2$ by means of triangles. All the results are obtained by applying 
our decoupled preconditioning technique. 
In each case we display the numerical rates of convergence versus the degrees of
freedom $DOF$. We observe  that, as expected, in the case of the MINI--BDM(1) and the 
$\mathbb P_2\mbox{-iso-}\mathbb P_1$--BDM(1) couplings, the convergence is linear  
for the velocities in both the Stokes and the Darcy domains. The Taylor-Hood--RT(1) scheme 
provides a quadratic convergence for the Stokes and Darcy velocity unknowns. 
We fixed the tolerance parameter for the outer MINRES method at $10^{-6}$ and 
checked empirically that the largest inner MINRES tolerance parameter  
that provides a convergence in agreement with the rates predicted by the theory is $10^{-2}$. 
All the results displayed here are obtained with this combination of tolerance parameters. 

\subsection{Performance of the iterative method}
In the following, we will denote  by $\mathbf{A}^h_{\mathrm{S}}$, $\mathbf{B}^h_{\mathrm{S}}$, 
$\mathbf{C}^h_{\mathrm{S}}$  and $\mathbf{M}^h_{\mathrm{S}}$ the matrix 
realizations of $A_{\mathrm{S}}^h$, $B_{\mathrm{S}}^h$, $C_{\mathrm{S}}^h$, 
and $I^h_{\mathrm{S}}$ respectively. Similarly,  $\mathbf{A}^h_{\mathrm{D}}$, $\mathbf{B}^h_{\mathrm{D}}$, 
$\mathbf{D}^h_{\mathrm{D}}$  and $\mathbf{M}^h_{\mathrm{D}}$ are the matrix 
realizations of $A_{\mathrm{D}}^h$, $B_{\mathrm{D}}^h$, $D_{\mathrm{D}}^h$, 
and $I^h_{\mathrm{D}}$ respectively. 

The numerical results were obtained using Matlab's own MINRES routine. For all experiments, 
the convergence is attained when the Euclidean norm of the relative residual is reduced by $10^{-6}$ for the outer MINRES 
while (as indicated above) the tolerance for the inner MINRES method is set to $10^{-2}$. The outer MINRES is applied to 
a linear system of equations with matrix
\[
\begin{pmatrix}
  \mathbf{A}^h_{\mathrm{S}} + \mathbf{C}^h_{\mathrm{S}} & (\mathbf{B}^h_{\mathrm{S}})^{\mathtt t}\\
  \mathbf{B}^h_{\mathrm{S}} & \mathbf{0}
 \end{pmatrix}.
\]
It is initialized with the solution 
of the Stokes problem with a non slip boundary condition $\Gamma_{\mathrm{S}}$ 
and an homogeneous Neumann boundary condition on $\Sigma$. The MINRES algorithm is  accelerated  
with one of the following preconditioners:
\[
\mathbf{\mathcal{P}}_{\mathrm{S}}^{\smallsetminus} :=\begin{pmatrix}
  (\mathbf{A}^h_{\mathrm{S}})^{-1}_{\smallsetminus} & \boldsymbol{0}\\
  \boldsymbol{0} & (\mathbf{M}^h_{\mathrm S})^{-1}
 \end{pmatrix},
\qquad 
\mathbf{\mathcal{P}}_{\mathrm{S}}^{\texttt{BPX}} :=\begin{pmatrix}
  (\mathbf{A}^h_{\mathrm{S}})^{-1}_{\texttt{BPX}} & \boldsymbol{0}\\
  \boldsymbol{0} & (\mathbf{M}^h_{\mathrm S})^{-1}
 \end{pmatrix},
\]
where the notation $(\mathbf{A}^h_{\mathrm{S}})^{-1}_{\smallsetminus}$ means 
that the linear systems of equations with matrix $\mathbf{A}^h_{\mathrm{S}}$ are solved by a direct solver (with the Matlab 
backslash command) while $(\mathbf{A}^h_{\mathrm{S}})^{-1}_{\texttt{BPX}}$ means that we use the 
Bramble-Pasciak-Xu \cite{BPX90, Xu92} preconditioner corresponding to the SPD vectorial 
Laplace matrix $\mathbf{A}^h_{\mathrm{S}}$.  

On the other hand, each application of $\mathbf{C}^h_{\mathrm{S}}$ to a vector requires the solution of 
a saddle point problem with matrix
\[
\begin{pmatrix}
  \mathbf{A}^h_{\mathrm{D}}  & (\mathbf{B}^h_{\mathrm{D}})^{\mathtt t}\\
  \mathbf{B}^h_{\mathrm{D}} & \mathbf{0}
 \end{pmatrix}.
\]
We again accomplish this task applying  the MINRES method preconditioned with one of the following symmetric and 
block diagonal matrices: 
\[
\mathbf{\mathcal{P}}_{\mathrm{D}}^{0} :=
 \begin{pmatrix}
  (\mathbf{A}^h_{\mathrm{D}} + \mathbf{D}^h_{\mathrm{D}})^{-1}_{\smallsetminus} & \mathbf{0}\\
  \mathbf{0} & (\mathbf{M}^h_\mathrm{D})^{-1}
 \end{pmatrix},
\quad 
\mathbf{\mathcal{P}}_{\mathrm{D}}^{\smallsetminus} := \begin{pmatrix}
  (\mathbf{P}^h_\mathrm{D})^{-1}_{\smallsetminus} & \mathbf{0}\\
  \mathbf{0} & (\mathbf{M}^h_\mathrm{D})^{-1}
 \end{pmatrix},
\]
\[
\mathbf{\mathcal{P}}_{\mathrm{D}}^{\texttt{BPX}}:=
 \begin{pmatrix}
  (\mathbf{P}^h_\mathrm{D})^{-1}_{\texttt{BPX}} & \mathbf{0}\\
  \mathbf{0} & (\mathbf{M}^h_\mathrm{D})^{-1}
 \end{pmatrix}.
\]
In the definition of the preconditioner $\mathbf{\mathcal{P}}_{\mathrm{D}}^{0}$, 
$(\mathbf{A}^h_{\mathrm{D}} + \mathbf{D}^h_{\mathrm{D}})^{-1}_{\smallsetminus}$ means that we simply use 
a direct solver for  
$\mathbf{A}^h_{\mathrm{D}} + \mathbf{D}^h_{\mathrm{D}}$ with the aid of  
the backslash Matlab command. The preconditioners 
$\mathbf{\mathcal{P}}_{\mathrm{D}}^{\smallsetminus}$ and $\mathbf{\mathcal{P}}_{\mathrm{D}}^{\texttt{BPX}}$ are 
obtained by substituting $(\mathbf{A}^h_{\mathrm{D}} + \mathbf{D}^h_{\mathrm{D}})^{-1}$ 
in $\mathbf{\mathcal{P}}_{\mathrm{D}}^{0}$ by the Hiptmair and Xu preconditioner 
$(\mathbf{P}^h_\mathrm{D})^{-1}$. The subscript $\smallsetminus$ in 
$(\mathbf{P}_\mathrm{D})^{-1}_{\smallsetminus}$ means that we solve the underlying 
Laplace problems with a direct solver, with the Matlab backslash command,  
and $(\mathbf{P}^h_\mathrm{D})^{-1}_{\texttt{BPX}}$ means that we use the well-known 
BPX-preconditioner \cite{BPX90, Xu92} for $(\mathbf{L}_h)^{-1}$ and $( -\Delta_h )^{-1}$.

In the cases where the mass matrix is diagonal the action of its inverse can be 
explicitly computed. In the other cases, one simple and effective strategy consists in substituting the action 
of the inverse of the mass matrix by one sweep of the  symmetric Gauss-Seidel method.

\begin{table}[ht!]
\footnotesize
\begin{center}
\begin{tabular}{r |c ||c |c |c |c |c |c }

$DOF$ & $h$ & $\mathbf{ \mathcal{P} }_{ \mathrm{S} }^{ \smallsetminus }(\mathbf{ \mathcal{P} }_{ \mathrm{D} }^{0})$ 
 & $\mathbf{ \mathcal{P} }_{ \mathrm{S} }^{ \smallsetminus }(\mathbf{ \mathcal{P} }_{ \mathrm{D} }^{ \smallsetminus })$ 
 & $\mathbf{ \mathcal{P} }_{ \mathrm{S} }^{ \smallsetminus } (\mathbf{\mathcal{P}}_{\mathrm{D}}^{ \mathtt{BPX} })$ 
 & $\mathbf{ \mathcal{P} }_{ \mathrm{S} }^{ \mathtt{BPX} }(\mathbf{\mathcal{P}}_{ \mathrm{D} }^{0} )$ 
 & $\mathbf{ \mathcal{P} }_{ \mathrm{S} }^{ \mathtt{BPX} }(\mathbf{\mathcal{P}}_{ \mathrm{D} }^{\smallsetminus})$ 
 & $\mathbf{ \mathcal{P} }_{ \mathrm{S} }^{ \mathtt{BPX} }(\mathbf{\mathcal{P}}_{ \mathrm{D} }^{\mathtt{BPX}} )$ \cr

\hline

543 & 1/8 & 26(4) & 26(26) & 26(29) & 56(4) & 56(26) & 56(29) \cr

2043 & 1/16 & 32(4) & 32(30) & 32(46) & 84(4) & 84(30) & 84(46) \cr

7923 & 1/32 & 40(4) & 40(33) & 40(62) & 121(4) & 121(33) & 121(62) \cr

31203 & 1/64 & 46(4) & 46(38) & 46(77) & 144(4) & 144(38) & 144(77) \cr

123843 & 1/128 & 50(4) & 50(42) & 50(91) & 158(4) & 158(42) & 158(91) \cr

\hline

\end{tabular}
\caption{Number of iterations: MINI--BDM(1)}
\label{table:6}
\end{center}
\end{table}

\begin{table}[ht!]
\footnotesize
\begin{center}
\begin{tabular}{r |c ||c |c |c |c |c |c }

$DOF$ & $h$ & $\mathbf{ \mathcal{P} }_{ \mathrm{S} }^{ \smallsetminus }(\mathbf{ \mathcal{P} }_{ \mathrm{D} }^{0})$ 
 & $\mathbf{ \mathcal{P} }_{ \mathrm{S} }^{ \smallsetminus }(\mathbf{ \mathcal{P} }_{ \mathrm{D} }^{ \smallsetminus })$ 
 & $\mathbf{ \mathcal{P} }_{ \mathrm{S} }^{ \smallsetminus } (\mathbf{\mathcal{P}}_{\mathrm{D}}^{ \mathtt{BPX} })$ 
 & $\mathbf{ \mathcal{P} }_{ \mathrm{S} }^{ \mathtt{BPX} }(\mathbf{\mathcal{P}}_{ \mathrm{D} }^{0} )$ 
 & $\mathbf{ \mathcal{P} }_{ \mathrm{S} }^{ \mathtt{BPX} }(\mathbf{\mathcal{P}}_{ \mathrm{D} }^{\smallsetminus})$ 
 & $\mathbf{ \mathcal{P} }_{ \mathrm{S} }^{ \mathtt{BPX} }(\mathbf{\mathcal{P}}_{ \mathrm{D} }^{\mathtt{BPX}} )$ \cr

\hline

385 & 1/8 & 24(4) & 24(26) & 24(29) & 50(4) & 50(26) & 50(29) \cr

1423 & 1/16 & 30(4) & 30(29) & 30(46) & 80(4) & 80(29) & 80(46) \cr

5467 & 1/32 & 36(4) & 36(34) & 36(62) & 107(4) & 107(34) & 107(62) \cr

21427 & 1/64 & 42(4) & 42(38) & 42(76) & 130(4) & 130(38) & 130(76) \cr

84835 & 1/128 & 44(4) & 44(41) & 44(91) & 146(4) & 146(41) & 146(91) \cr

\hline

\end{tabular}
\caption{Number of iterations: $\mathbb{P}$2-iso-$\mathbb{P}$1--BDM(1)}
\label{table:7}
\end{center}
\end{table}

%
%
%
%
%
%
%

\begin{table}[ht!]
\footnotesize
\begin{center}
\begin{tabular}{r |c ||c |c }

$DOF$ & $h$ & $\mathbf{ \mathcal{P} }_{ \mathrm{S} }^{ \smallsetminus }(\mathbf{ \mathcal{P} }_{ \mathrm{D} }^{0})$ 
& $\mathbf{ \mathcal{P} }_{ \mathrm{S} }^{ \smallsetminus }(\mathbf{ \mathcal{P} }_{ \mathrm{D} }^{ \smallsetminus })$ \cr

\hline

887 & 1/8 & 28(5) & 28(28) \cr

3371 & 1/16 & 34(5) & 34(32) \cr

13139 & 1/32 & 38(5) & 38(36) \cr

31203 & 1/64 & 42(5) & 42(40) \cr

206147 & 1/128 & 44(5) & 44(44) \cr

\hline

\end{tabular}
\caption{Number of iterations: Taylor-Hood--RT(1)}
\label{table:8}
\end{center}
\end{table}

In tables \ref{table:6},  \ref{table:7} and \ref{table:8}, we list 
the number of iterations of the two nested MINRES methods with different combinations of 
preconditioners. The preconditioner in brackets is the one used for the inner MINRES. 
We show the number of outer MINRES iterations and the number in brackets 
is an average of the number of inner MINRES iterations. We observe that for different mesh sizes, the iterative method results 
in a uniform number of MINRES iterations. Therefore, the preconditioners  
are robust with respect to the mesh size, which agrees with the theoretical results. 

%

%

\end{document}